\documentclass[a4paper,12pt]{article}
\usepackage{amsmath}
\usepackage{amssymb}
\usepackage{amsthm}
\usepackage[latin1]{inputenc}     

\newcommand{\beq}{\begin{equation} }
\newcommand{\eqq}{\end{equation} }
\newcommand{\cuad}{{\sqcap\kern-.68em\sqcup}}

\newtheorem{remark}{Remark}[section]
\newcommand{\bremark}{\begin{remark} \em}
\newcommand{\eremark}{\end{remark} }

\def\beeq{\begin{equation}}
\def\eeq{\end{equation}}
\newcommand{\begeqaet}{\begin{eqnarray*}}
\newcommand{\eneqaet}{\end{eqnarray*}}

\usepackage{mathtools} 

\let\Section=\section
\def\section{\setcounter{equation}{0}\Section}
\newtheorem{Lem}{Lemma}[section]
\newtheorem{Thm}{Theorem}[section]
\newtheorem{Def}{Definition}[section]
\newtheorem{Prop}{Proposition}[section]
\newtheorem{Remark}{Remark}[section]

\begin{document}
\begin{center}{\bf\Large Boundary value problem with fractional p-Laplacian operator}\medskip

\bigskip

\bigskip

{C\'esar Torres}

Departamento de Matem\'aticas\\
Universidad Nacional de Trujillo\\
Av. Juan Pablo Segundo s/n, Trujillo-Per\'u

 {\sl  (ctl\_576@yahoo.es, ctorres@dim.uchile.cl)}

\end{center}

\medskip

\medskip
\medskip
\medskip
\medskip

\begin{abstract}
The aim of this paper is to obtain the existence of solution for the  fractional p-Laplacian Dirichlet problem with mixed derivatives
\begin{eqnarray*}
&{_{t}}D_{T}^{\alpha}\left(|_{0}D_{t}^{\alpha}u(t))|^{p-2}{_{0}}D_{t}^{\alpha}u(t)\right)  =  f(t,u(t)), \;t\in [0,T],\\
&u(0)  =  u(T) = 0,
\end{eqnarray*}
where $\frac{1}{p} < \alpha <1$, $1<p<\infty$ and $f:[0,T]\times \mathbb{R} \to \mathbb{R}$ is a Carath\'eodory function wich satisfies some growth conditions. We obtain the existence of nontrivial solution by using the Mountain Pass Theorem.


{\bf MSC}
\end{abstract}
\date{}

\setcounter{equation}{0}
\section{ Introduction}

Recently, a great attention has been focused on the study of boundary value problems (BVP) for fractional differential equations. They appear in mathematical models in different branches in Science as physics, chemistry, biology, geology, as well as, control theory, signal theory, nanoscience and so on \cite{DBZGJM, AKHSJT, IP, JSOAJTM, SSAKOM, YZ} and references therein. 

Physical models containing left and right fractional differential operators have recently renewed attention from scientists which is mainly due to applications as models for physical phenomena exhibiting anomalous diffusion. Specifically, the models involving a fractional differential oscillator equation, which contains a composition of left and right fractional derivatives, are proposed for the description of the processes of emptying the silo \cite{SLTB} and the heat flow through a bulkhead filled with granular material \cite{ES}, respectively. Their studies show that the proposed models based on fractional calculus are efficient and describe well the processes.

The existence and multiplicity of solutions for BVP for nonlinear fractional differential equations is extensively studied using various tools of nonlinear analysis as fixed point theorems, degree theory and the method of upper and lower solutions \cite{MBJNRR, MBACDS}. Very recently, it should be noted that critical point theory and variational methods have also turned out to be very effective tools in determining the existence of solutions of BVP for fractional differential equations. The idea behind them is trying to find solutions of a given boundary value problem by looking for critical points of a suitable energy functional defined on an appropriate function space. In the last 30 years, the critical point theory has become a wonderful tool in studying the existence of solutions to differential equations with variational structures, we refer the reader to the books due to Mawhin and Willem \cite{JMMW}, Rabinowitz \cite{PR}, Schechter \cite{MS} and papers \cite{VEJR, FJYZ0, FJYZ, CT, CT1, CT2, WXJXZL, YZ}.

The p-Laplacian operator was considered in several recent works. It arises in the modelling of different physical and natural phenomena; non-Newtonian mechanics, nonlinear elasticity and glaciology, combustion theory, population biology, nonlinear flow laws, system of Monge-Kantorovich partial differential equations. There exists a very large number of papers devoted to the existence of solutions of the p-Laplacian operator in which the authors used bifurcation, variational methods, sub-super solutions, degree theory, in order to prove the existence of solutions of this nonlinear operator, for detail see \cite{GDPJJM}.
 
Motivated by these previous works, we consider the solvability of the Dirichlet problem with mixed fractional derivatives 
\begin{eqnarray}\label{I01}
&{_{t}}D_{T}^{\alpha}\left(|_{0}D_{t}^{\alpha}u(t))|^{p-2}{_{0}}D_{t}^{\alpha}u(t)\right)  =  f(t,u(t)), \;t\in [0,T],\nonumber\\
&u(0)  =  u(T) = 0,
\end{eqnarray}
where $1<p<\infty$, $\frac{1}{p}< \alpha < 1$ and we assume that $f: [0,T]\times \mathbb{R} \to \mathbb{R}$ is a Carath\'eodory function satisfying:
\begin{itemize}
\item[\fbox{$f_1$}] There exists $C>0$ and $1<q<\infty$, such that
$$
|f(t,\xi)| \leq C(1 + |\xi|^{q-1})\;\;\mbox{such that for a.e.}\;t\in [0,T],\;\xi \in \mathbb{R}
$$
\item[\fbox{$f_2$}] There exists $\mu >p$ and $r>0$ such that for a.e. $t\in [0,T]$ and $r\in \mathbb{R}$, $|\xi| \geq r$
$$
0< \mu F(t,\xi) \leq \xi f(t,\xi),
$$
where $F(t,\xi) = \int_{0}^{\xi}f(t,\sigma)d\sigma$.
\item[\fbox{$f_3$}] $\lim_{\xi \to 0} \frac{f(t,\xi)}{|\xi|^{p-1}} = 0$ uniformly for a.e. $t\in [0,T]$.
\end{itemize}

We say that $u\in E_{0}^{\alpha ,p}$ is a weak solution of problem (\ref{I01}), if
$$
\int_{0}^{T} |{_{0}}D_{t}^{\alpha}u(t)|^{p-2}{_{0}}D_{t}^{\alpha}u(t) {_{0}}D_{t}^{\alpha} \varphi (t)dt = \int_{0}^{T} f(t,u(t))\varphi (t)dt,
$$
for any $\varphi \in E_{0}^{\alpha,p}$, where space $E_{0}^{\alpha ,p}$ will be introduced in Section $\S$ 2.

Let $I: E_{0}^{\alpha, p} \to \mathbb{R}$ the functional associated to (\ref{I01}), defined by
\begin{equation}\label{I02}
I(u) = \frac{1}{p}\int_{0}^{T} |{_{0}}D_{t}^{\alpha}u(t)|^{p}dt - \int_{\mathbb{R}}F(t,u(t))dt
\end{equation}
under our assumption $I\in C^1$ and we have
\begin{equation}\label{I03}
I'(u)v = \int_{0}^{T} |{_{0}}D_{t}^{\alpha}u(t)|^{p-2}{_{0}}D_{t}^{\alpha}u(t){_{0}}D_{t}^{\alpha}v(t)dt - \int_{0}^{T}f(t,u(t))v(t)dt.
\end{equation}
Moreover critical points of $I$ are weak solutions of problem (\ref{I01}).

Using the Mountain pass Theorem, we get our main result.
\begin{Thm}\label{main}
Suppose that $f$ satisfies $(f_1) - (f_3)$. If $p<q< \infty$ then the problem (\ref{I01}) has a nontrivial weak solution in $E_{0}^{\alpha ,p}$. 
\end{Thm}

The rest of the paper is organized as follows: In Section \S 2 we present preliminaries on
fractional calculus and we introduce the functional setting of the problem. In Section \S 3 we prove Theorem \ref{main}.

\section{Fractional Calculus}

In this section we introduce some basic definitions of fractional calculus which are used further in this paper.  For the proof see \cite{AKHSJT, IP, SSAKOM}.


Let $u$ be a function defined on $[a,b]$. The left (right ) Riemann-Liouville fractional integral of order $\alpha >0$ for function $u$ is defined by
\begin{eqnarray*}
&_{a}I_{t}^{\alpha}u(t) = \frac{1}{\Gamma (\alpha)}\int_{a}^{t} (t-s)^{\alpha - 1}u(s)ds,\;t\in [a,b],\\
&_{t}I_{b}^{\alpha}u(t) = \frac{1}{\Gamma (\alpha)}\int_{t}^{b}(s-t)^{\alpha -1}u(s)ds,\;t\in [a,b],
\end{eqnarray*}
provided in both cases that the right-hand side is pointwise defined on $[a,b]$.


The left and right Riemann - Liouville fractional derivatives of order $\alpha >0$ for function $u$ denoted by $_{a}D_{t}^{\alpha}u(t)$ and $_{t}D_{b}^{\alpha}u(t)$, respectively, are defined by
\begin{eqnarray*}
&_{a}D_{t}^{\alpha}u(t) = \frac{d^{n}}{dt^{n}}{_{a}}I_{t}^{n-\alpha}u(t),\\
&_{t}D_{b}^{\alpha}u(t) = (-1)^{n}\frac{d^{n}}{dt^{n}}{ _{t}}I_{b}^{n-\alpha}u(t),
\end{eqnarray*}
where $t\in [a,b]$, $n-1 \leq \alpha < n$ and $n\in \mathbb{N}$.

The left and right Caputo fractional derivatives are defined via the above Riemann-Liouville fractional derivatives \cite{AKHSJT}. In particular, they are defined for the function belonging to the space of absolutely continuous function, namely, If $\alpha \in (n-1,n)$ and $u\in AC^{n}[a,b]$, then the left and right Caputo fractional derivative of order $\alpha$ for function $u$ denoted by $_{a}^{c}D_{t}^{\alpha}u(t)$ and $_{t}^{c}D_{b}^{\alpha}u(t)$ respectively, are defined by 
\begin{eqnarray*}
&& _{a}^{c}D_{t}^{\alpha}u(t) = {_{a}}I_{t}^{n-\alpha}u^{(n)}(t) = \frac{1}{\Gamma (n-\alpha)}\int_{a}^{t} (t-s)^{n-\alpha -1}u^{n}(s)ds,\\
&&_{t}^{c}D_{b}^{\alpha}u(t) = (-1)^{n} {_{t}}I_{b}^{n-\alpha}u^{(n)}(t) = \frac{(-1)^{n}}{\Gamma (n-\alpha)}\int_{t}^{b} (s-t)^{n-\alpha-1}u^{(n)}(s)ds.
\end{eqnarray*}

The Riemann-Liouville fractional derivative and the Caputo fractional derivative are connected with each other by the following relations
\begin{Thm}\label{RL-C}
Let $n \in \mathbb{N}$ and $n-1 < \alpha < n$. If $u$ is a function defined on $[a,b]$ for which the Caputo fractional derivatives $_{a}^{c}D_{t}^{\alpha}u(t)$ and $_{t}^{c}D_{b}^{\alpha}u(t)$ of order $\alpha$ exists together with the Riemann-Liouville fractional derivatives $_{a}D_{t}^{\alpha}u(t)$ and $_{t}D_{b}^{\alpha}u(t)$, then
\begin{eqnarray*}
_{a}^{c}D_{t}^{\alpha}u(t) & = & _{a}D_{t}^{\alpha}u(t) - \sum_{k=0}^{n-1} \frac{u^{(k)}(a)}{\Gamma (k-\alpha + 1)} (t-a)^{k-\alpha}, \quad t\in [a,b], \\
_{t}^{c}D_{b}^{\alpha}u(t) & = & _{t}D_{b}^{\alpha}u(t) - \sum_{k=0}^{n-1} \frac{u^{(k)}(b)}{\Gamma (k-\alpha + 1)} (b-t)^{k-\alpha},\quad t\in [a,b].
\end{eqnarray*}
In particular, when $0<\alpha < 1$, we have
\begin{equation}\label{RL-C01}
_{a}^{c}D_{t}^{\alpha}u(t) =  {_{a}}D_{t}^{\alpha}u(t) -  \frac{u(a)}{\Gamma (1-\alpha)} (t-a)^{-\alpha}, \quad t\in [a,b]
\end{equation}
and
\begin{equation}\label{RL-C02}
_{t}^{c}D_{b}^{\alpha}u(t) = {_{t}}D_{b}^{\alpha}u(t) - \frac{u(b)}{\Gamma (1-\alpha)}(b-t)^{-\alpha},\quad t\in [a,b].
\end{equation}

\end{Thm}

Now we consider some properties of the Riemann-Liouville fractional integral and derivative operators.

\begin{itemize}
\item[(1)] 
\begin{eqnarray*}
&&_{a}I_{t}^{\alpha}(_{a}I_{t}^{\beta}u(t)) =  {_{a}}I_{t}^{\alpha + \beta}u(t)\;\;\mbox{and}\\
&&_{t}I_{b}^{\alpha}(_{t}I_{b}^{\beta}u(t)) = { _{t}}I_{b}^{\alpha + \beta}u(t)\;\;\forall \alpha, \beta >0,
\end{eqnarray*}

\item[(2)] {\bf Left inverse.} Let $u \in L^{1}[a,b]$ and $\alpha >0$,
\begin{eqnarray*}
&&_{a}D_{t}^{\alpha}(_{a}I_{t}^{\alpha}u(t)) =  u(t),\;\mbox{a.e.}\;t\in[a,b]\;\;\mbox{and}\\
&&_{t}D_{b}^{\alpha}(_{t}I_{b}^{\alpha}u(t)) =  u(t),\;\mbox{a.e.}\;t\in[a,b].
\end{eqnarray*}

\item[(3)] For $n-1\leq \alpha < n$, if the left and right Riemann-Liouville fractional derivatives $_{a}D_{t}^{\alpha}u(t)$ and $_{t}D_{b}^{\alpha}u(t)$, of the function $u$ are integral on $[a,b]$, then
\begin{eqnarray*}
_{a}I_{t}^{\alpha}(_{a}D_{t}^{\alpha}u(t)) & = & u(t) - \sum_{k = }^{n} [_{a}I_{t}^{k-\alpha}u(t)]_{t=a} \frac{(t-a)^{\alpha -k}}{\Gamma (\alpha - k + 1)},\\
_{t}I_{b}^{\alpha}(_{t}D_{b}^{\alpha}u(t)) & = & u(t) - \sum_{k=1}^{n}[_{t}I_{n}^{k-\alpha}u(t)]_{t=b}\frac{(-1)^{n-k}(b-t)^{\alpha - k}}{\Gamma (\alpha - k +1)},
\end{eqnarray*}
for $t\in [a,b]$.

\item[(4)] {\bf Integration by parts}
\begin{equation}\label{FCeq1}
\int_{a}^{b}[_{a}I_{t}^{\alpha}u(t)]v(t)dt = \int_{a}^{b}u(t)_{t}I_{b}^{\alpha}v(t)dt,\;\alpha >0,
\end{equation}
provided that $u\in L^{p}[a,b]$, $v\in L^{q}[a,b]$ and
$$
p\geq 1,\;q\geq 1\;\;\mbox{and}\;\;\frac{1}{p}+\frac{1}{q} < 1+\alpha \;\;\mbox{or}\;\; p \neq 1,\;q\neq 1\;\;\mbox{and}\;\;\frac{1}{p} + \frac{1}{q} = 1+\alpha.
$$

\begin{equation}\label{FCeq2}
\int_{a}^{b} [_{a}D_{t}^{\alpha}u(t)]v(t)dt = \int_{a}^{b}u(t)_{t}D_{b}^{\alpha}v(t)dt,\;\;0<\alpha \leq 1,
\end{equation}
provided the boundary conditions
\begin{eqnarray*}
&u(a) = u(b) = 0,\;u'\in L^{\infty}[a,b],\;v\in L^{1}[a,b]\;\;\mbox{or}\\
&v(a) = v(b) = 0,\;v' \in L^{\infty}[a,b], \;u \in L^{1}[a,b],
\end{eqnarray*}
are fulfilled.
\end{itemize}
\subsection{Fractional Derivative Space}

In order to establish a variational structure for BVP (\ref{I01}), it is necessary to construct appropriate function spaces. For this setting we take some results from \cite{FJYZ0, FJYZ, YZ}.

Let us recall that for any fixed $t\in [0,T]$ and $1\leq p <\infty$,
\begin{eqnarray*}
\|u\|_{L^{p}[0,t]}  =  \left( \int_{0}^{t} |u(s)|^{p}ds \right)^{1/p},\;
\|u\|_{L^{p}}  = \left( \int_{0}^{T} |u(s)|^{p}ds \right)^{1/p}\;\;\mbox{and}\;\;
\|u\|_{\infty} = \max_{t\in [0,T]}|u(t)|.
\end{eqnarray*}

\begin{Def}\label{FC-FEdef1}
Let $0< \alpha \leq 1$ and $1<p<\infty$. The fractional derivative spaces $E_{0}^{\alpha ,p}$ is defined by
\begin{eqnarray*}
E_{0}^{\alpha , p} &= & \{u\in L^{p}[0,T]/\;\;_{0}D_{t}^{\alpha}u \in L^{p}[0,T]\;\mbox{and}\;u(0) = u(T) = 0\}\\
&= & \overline{C_{0}^{\infty}[0,T]}^{\|.\|_{\alpha ,p}}.
\end{eqnarray*}
where $\|.\|_{\alpha ,p}$ is defined by
\begin{equation}\label{FC-FEeq1}
\|u\|_{\alpha ,p}^{p} = \int_{0}^{T} |u(t)|^{p}dt + \int_{0}^{T}|_{0}D_{t}^{\alpha}u(t)|^{p}dt.
\end{equation}
\end{Def}

\begin{Remark}\label{RL-Cnta}
For any $u\in E_{0}^{\alpha, p}$, nothing the fact that $u(0) = 0$, we have ${^{c}_{0}}D_{t}^{\alpha}u(t) = {_{0}}D_{t}^{\alpha}u(t)$, $t\in [0,T]$ according to (\ref{RL-C01}).
\end{Remark}
\begin{Prop}\label{FC-FEprop1}
\cite{FJYZ0} Let $0< \alpha \leq 1$ and $1 < p <\infty$. The fractional derivative space $E_{0}^{\alpha , p}$ is a reflexive and separable Banach space.
\end{Prop}

We recall some properties of the fractional space $E_{0}^{\alpha ,p}$. 
\begin{Lem}\label{FC-FElem1}
\cite{FJYZ0} Let $0< \alpha \leq 1$ and $1\leq p < \infty$. For any $u\in L^{p}[0,T]$ we have
\begin{equation}\label{FC-FEeq2}
\|_{0}I_{\xi}^{\alpha}u\|_{L^{p}[0,t]}\leq \frac{t^{\alpha}}{\Gamma(\alpha + 1)} \|u\|_{L^{p}[0,t]},\;\mbox{for}\;\xi\in [0,t],\;t\in[0,T].
\end{equation}
\end{Lem}


\begin{Prop}\label{FC-FEprop3}
\cite{FJYZ} Let $0< \alpha \leq 1$ and $1 < p < \infty$. For all $u\in E_{0}^{\alpha ,p}$,  we have
\begin{equation}\label{FC-FEeq3}
\|u\|_{L^{p}} \leq \frac{T^{\alpha}}{\Gamma (\alpha +1)} \|_{0}D_{t}^{\alpha}u\|_{L^{p}}.
\end{equation}
If $\alpha > 1/p$ and $\frac{1}{p} + \frac{1}{q} = 1$, then
\begin{equation}\label{FC-FEeq4}
\|u\|_{\infty} \leq \frac{T^{\alpha -1/p}}{\Gamma (\alpha)((\alpha - 1)q +1)^{1/q}}\|_{0}D_{t}^{\alpha}u\|_{L^{p}}.
\end{equation}
\end{Prop}

\begin{Remark}\label{embb}
Let $1/p< \alpha \leq 1$, if $u\in E_{0}^{\alpha, p}$, then $u\in L^{q}[0,T]$ for $q\in [p, +\infty]$. In fact
\begin{eqnarray*}
\int_{0}^{T} |u(t)|^{q}dt &=& \int_{0}^{T} |u(t)|^{q-p}|u(t)|^{p}dt\\
& \leq & \|u\|_{\infty}^{q-p} \|u\|_{L^{p}}^{p}.
\end{eqnarray*}
In particular the embedding $E_{0}^{\alpha ,p} \hookrightarrow L^{q}[0,T]$ is continuos for all $q\in [p, +\infty]$. 
\end{Remark}

\noindent 
According to (\ref{FC-FEeq3}), we can consider in $E_{0}^{\alpha ,p}$ the following norm
\begin{equation}\label{FC-FEeq5}
\|u\|_{\alpha ,p} = \|_{0}D_{t}^{\alpha}u\|_{L^{p}},
\end{equation}
and (\ref{FC-FEeq5}) is equivalent to (\ref{FC-FEeq1}).

\begin{Prop}\label{FC-FEprop4}
\cite{FJYZ} Let $0< \alpha \leq 1$ and $1 < p < \infty$. Assume that $\alpha > \frac{1}{p}$ and $\{u_{k}\} \rightharpoonup u$ in $E_{0}^{\alpha ,p}$. Then $u_{k} \to u$ in $C[0,T]$, i.e.
$$
\|u_{k} - u\|_{\infty} \to 0,\;k\to \infty.
$$
\end{Prop}

Now, we are going to prove that $E_{0}^{\alpha , p}$ is uniformly convex, for this fact we consider the following tools (see \cite{RAJF} for more details).
\begin{itemize}
\item[(1)] {\bf Reverse H\"older Inequality:} Let $0<p<1$, so that $p' = \frac{p}{p-1} <0$. If $u\in L^{p}(\Omega)$ and
$$
0< \int_{\Omega}|g(x)|^{p'}dx < \infty,
$$
then
\begin{equation}\label{HI}
\int_{\Omega} |f(x)g(x)|dx \geq \left( \int_{\Omega} |f(x)|^{p}dx \right)^{1/p} \left(\int_{\Omega} |g(x)|^{p'} \right)^{1/p'}.
\end{equation}
\item[(2)] {\bf Reverse Minkowski inequality:} Let $0<p<1$. If $u,v\in L^{p}(\Omega)$, the
\begin{equation}\label{MI}
\||u| + |v|\|_{L^p} \geq \|u\|_{p} + \|v\|_{p}
\end{equation}
\item[(3)] Let $z,w\in \mathbb{C}$. If $1< p \leq 2$ and $p'=\frac{p}{p-1}$, then
\begin{equation}\label{p1}
\left|\frac{z+w}{2} \right|^{p'} + \left| \frac{z-w}{2} \right|^{p'} \leq \left(\frac{1}{2}|z|^p + \frac{1}{2}|w|^p  \right)^{1/(p-1)}.
\end{equation}
If $2\leq p < \infty$, then
\begin{equation}\label{p2}
\left| \frac{z+w}{2} \right|^{p} + \left| \frac{z-w}{2} \right|^{p} \leq \frac{1}{2}|z|^p + \frac{1}{2}|w|^p
\end{equation}
\end{itemize}

\begin{Lem}\label{lemu}
$(E_{0}^{\alpha ,p}, \|.\|_{\alpha , p})$ is uniformly convex.
\end{Lem}

\noindent 
{\bf Proof.} Let $u,v \in E_{0}^{\alpha ,p}$ satisfy $\|u\|_{\alpha ,p} = \|v\|_{\alpha ,p} = 1$ and $\|u-v\|_{\alpha ,p} \geq \epsilon$, where $\epsilon \in (0,2)$.

\noindent
{\bf Case $p\geq 2$.} By (\ref{p2}), we have
\begin{eqnarray}\label{U01}
\left\| \frac{u+v}{2} \right\|_{\alpha ,p}^{p} + \left\| \frac{u-v}{2} \right\|_{\alpha ,p}^{p}\!\!\!\!\!& = &\!\!\! \int_{0}^{T} \left| \frac{{_{0}}D_{t}^{\alpha}u(t) + {_{0}}D_{t}^{\alpha}v(t)}{2}  \right|^pdt + \int_{0}^{T} \left| \frac{{_{0}}D_{t}^{\alpha}u(t) - {_{0}}D_{t}^{\alpha}v(t)}{2}  \right|^{p}dt\nonumber\\
&\leq& \frac{1}{2} \int_{0}^{t} |{_{0}}D_{t}^{\alpha}u(t)|^pdt + \frac{1}{2}\int_{0}^{T} |{_{0}}D_{t}^{\alpha}v(t)|^{p}dt\nonumber\\
& = & \frac{1}{2}\|u\|_{\alpha ,p}^{p} + \frac{1}{2}\|v\|_{\alpha ,p}^{p} = 1.
\end{eqnarray}
It follows from (\ref{U01}) that
$$
\left\| \frac{u+v}{2} \right\|_{\alpha ,p}^{p}  \leq 1 - \frac{\epsilon^p}{2^p}.
$$
Taking $\delta = \delta (\epsilon)$ such that $1-(\epsilon/2)^2 = (1-\delta)^{p}$, we obtain that
$$
\left\| \frac{u+v}{2}  \right\|_{\alpha ,p} \leq (1-\delta).
$$

\noindent
{\bf Case $1<p<2$.} First, note that
$$
\|u\|_{\alpha ,p}^{p'} = \left( \int_{0}^{T} \left( |{_{0}}D_{t}^{\alpha}u(t)|^{p'}  \right)^{p-1}dt \right)^{\frac{1}{p-1}},
$$
where $p' = \frac{p}{p-1}$. Using the reverse Minkowski inequality (\ref{MI}) and the inequality (\ref{p1}), we get
\begin{eqnarray}\label{U02}
&&\left\| \frac{u+v}{2} \right\|_{\alpha ,p}^{p'} + \left\| \frac{u-v}{2} \right\|_{\alpha ,p}^{p'}  \nonumber \\
&& =  \left[ \int_{0}^{T} \left( \left| \frac{{_{0}}D_{t}^{\alpha}u(t) + {_{0}}D_{t}^{\alpha}v(t) }{2} \right|^{p'}  \right)^{p-1} dt\right]^{\frac{1}{p-1}}  +  \left[ \int_{0}^{T} \left( \left| \frac{{_{0}}D_{t}^{\alpha}u(t) - {_{0}}D_{t}^{\alpha}v(t) }{2} \right|^{p'}  \right)^{p-1}dt \right]^{\frac{1}{p-1}}\nonumber\\
&&\leq \left[ \int_{0}^{T} \left( \left| \frac{{_{0}}D_{t}^{\alpha}u(t) + {_{0}}D_{t}^{\alpha}v(t)}{2}  \right|^{p'} + \left| \frac{{_{0}}D_{t}^{\alpha}u(t) - {_{0}}D_{t}^{\alpha}v(t)}{2}  \right|^{p'}  \right)^{p-1}dt  \right]^{\frac{1}{p-1}}\nonumber\\
&& \leq \left[ \int_{0}^{T} \left( \frac{|{_{0}}D_{t}^{\alpha}u(t)|^p}{2} + \frac{|{_{0}}D_{t}^{\alpha}v(t)|^p}{2}\right)dt \right]^{p' -1}\nonumber\\
&& = \left(\frac{1}{2}\|u\|_{\alpha ,p}^{p} + \frac{1}{2}\|v\|_{\alpha ,p}^{p} \right)^{p' -1} = 1.
\end{eqnarray} 
By (\ref{U02}), we have
$$
\left\| \frac{u+v}{2} \right\|_{\alpha ,p}^{p'} \leq 1 - \frac{\epsilon^{p'}}{2^{p'}}.
$$
Taking $\delta = \delta(\epsilon)$ such that $1-(\epsilon /2)^{p'} = (1-\delta)^{p'}$, we get the desired claim. $\Box$

\section{Proof of Theorem \ref{main}}

Through this section we consider: $p<q$ and $\frac{1}{p} < \alpha \leq 1$. For $u\in E_{0}^{\alpha ,p}$ we define
$$
J(u) = \frac{1}{p} \int_{0}^{T} |{_{0}}D_{t}^{\alpha}u(t)|^{p}dt,\;\;H(u) = \int_{0}^{T}F(t,u(t))dt,
$$
and
$$
I(u) = J(u) - H(u).
$$
Obviously, the energy functional $I: E_{0}^{\alpha,p} \to \mathbb{R}$ associated with problem (\ref{I01}) is well defined.
\begin{Lem}\label{MR1lem}
If $f$ satisfies assumption $(f_1)$, then the functional $H\in C^1(E_{0}^{\alpha,p}, \mathbb{R})$ and
$$
\langle H'(u), v \rangle = \int_{0}^{T} f(t,u(t))v(t)dt\;\;\mbox{for all}\;\;u,v\in E_{0}^{\alpha,p}.
$$
\end{Lem}

\noindent 
{\bf Proof.}
\begin{itemize}
\item[(i)] $H$ is G\^ateaux-differentiable in $E_{0}^{\alpha ,p}$.

Let $u,v\in E_{0}^{\alpha ,p}$. For each $t\in [0,T]$ and $0< |\sigma| <1$, by the mean value theorem, there exists $0< \delta <1$,
\begin{eqnarray*}
\frac{1}{\sigma} (F(t,u+\sigma v) - F(t,u)) & = & \frac{1}{\sigma} \int_{0}^{u+\sigma v} f(t,s)ds - \frac{1}{\sigma} \int_{0}^{u}f(t,s)ds\\
& = & \frac{1}{\sigma} \int_{u}^{u + \sigma v} f(t,s)ds = f(t, u+\delta \sigma v)v.
\end{eqnarray*}
By ($f_1$) and Young's inequality, we get
\begin{eqnarray*}
|f(t,u+\delta \sigma v) v| &\leq& C(|v| + | u + \delta \sigma v|^{q-1}|v| )\\
&\leq& C(2|v|^q + |u+\delta \sigma v|^q + 1) \\
&\leq& a2^q (|v|^q + |u|^q + 1).
\end{eqnarray*}
Since $q>1 $, by  (\ref{FC-FEeq3}) we have $u,v \in L^{q}[0,T]$. Moreover, the Lebesgue Dominated Convergence Theorem implies
\begin{eqnarray*}
\lim_{\sigma \to 0} \frac{1}{\sigma} (H(u+\sigma v) - H(u)) & = & \lim_{\sigma \to 0} \int_{0}^{T} f(t, u + \delta \sigma v)vdt\\
& = & \int_{0}^{T} \lim_{\sigma \to 0}f(t, u + \delta \sigma v)vdt = \int_{0}^{T}f(t,u)vdt.
\end{eqnarray*}

\item[(ii)] Continuity of G\^ateaux-derivative.

Let $\{u_n\}, u\in E_{0}^{\alpha ,p}$ such that $u_n \to u$ strongly in $E_{0}^{\alpha ,p}$ as $n\to \infty$. Without loss of generality, we assume that $u_n(t) \to u(t)$ a.e. in $[0,T]$. By ($f_1$), for any $I \subset [0,T]$,
\begin{eqnarray}\label{MR01}
\int_{I} |f(t,u_n)|^{q'}dt &\leq& C^{q'}\int_{I} (1 + |u_n|^{q-1})^{q'}dt\nonumber\\
&\leq& C^{q'}2^{q'}\int_{I} (1 + |u_n|^q)dt\nonumber\\
&\leq & \overline{C}[\mu(I) + \|u_n\|_{\infty}^q \mu(I)],
\end{eqnarray}
where $\mu$ denotes the Lebesgue measure of $I$. It follows from (\ref{MR01}) that the sequence $\{|f(t,u_n) - f(t,u)|^{q'}\}$ is uniformly bounded and equi-integrable in $L^1[0,T]$. The Vitali Convergence Theorem implies
$$
\lim_{n\to \infty} \int_{0}^{T} |f(t,u_n) - f(t,u)|^{q'}dt = 0.
$$
Thus, by H\"older inequality and Remark \ref{embb}, we obtain
\begin{eqnarray*}
\|H'(u_n) - H(u)\|_{(E_{0}^{\alpha ,p})^{*}} & = & \sup_{v\in E_{0}^{\alpha,p}, \|v\|_{\alpha,p}=1} \left| \int_{0}^{T} (f(t,u_n) - f(t,u))vdt  \right|\\
&\leq & \|f(t,u_n) - f(t,u)\|_{L^{q'}} \|v\|_{L^q}\\
&\leq& K\|f(t,u_n) - f(t,u)\|_{L^{q'}}\\
&\to & 0,
\end{eqnarray*}
as $n\to \infty$. Hence, we complete the proof of Lemma. $\Box$
\end{itemize}

\begin{Lem}\label{MR2lem}
The functional $J \in C^{1}(E_{0}^{\alpha ,p}, \mathbb{R})$ and
$$
\langle J'(u), v \rangle = \int_{0}^{T} |{_{0}}D_{t}^{\alpha}u(t)|^{p-2}{_{0}}D_{t}^{\alpha}u(t){_{0}}D_{t}^{\alpha}v(t)dt,
$$
for all $u,v\in E_{0}^{\alpha ,p}$. Moreover, for each $u\in E_{0}^{\alpha ,p}$, $J'(u) \in (E_{0}^{\alpha,p})^{*}$, where $(E_{0}^{\alpha ,p})^{*}$ denotes the dual of $E_{0}^{\alpha,p}$. 
\end{Lem}

\noindent
{\bf Proof.}

First, it is easy to see that
\begin{equation}\label{MR02}
\langle J'(u), v \rangle = \int_{0}^{T} |{_{0}}D_{t}^{\alpha}u(t)|^{p-2}{_{0}}D_{t}^{\alpha}u(t) {_{0}}D_{t}^{\alpha}v(t)dt,
\end{equation}
for all $u,v \in E_{0}^{\alpha ,p}$. It follows from (\ref{MR02}) that for each $u\in E_{0}^{\alpha ,p}$, $J'(u) \in (E_{0}^{\alpha ,p})^{*}$.

Next, we prove that $J \in C^1(E_{0}^{\alpha ,p}, \mathbb{R})$. For the proof we need the following inequalities, (see \cite{GDPJJM})  
\begin{itemize}
\item[(i)] If $p\in [2,\infty)$ then it holds
\begin{equation}\label{MR03}
\left| |z|^{p-2}z - |y|^{p-2}y  \right| \leq \beta|z-y|(|z| + |y|)^{p-2}\;\;\mbox{for all}\;y,z\in \mathbb{R},
\end{equation}
with $\beta$ independent of $y$ and $z$;
\item[(ii)] If $p\in (1,2]$ then it holds:
\begin{equation}\label{MR04}
\left| |z|^{p-2}z - |y|^{p-2}y \right| \leq \beta |z-y|^{p-1}\;\;\mbox{for all}\;y,z\in \mathbb{R},
\end{equation} 
with $\beta$ independent of $y$ and $z$.
\end{itemize}

We define $g:E_{0}^{\alpha ,p} \to L^{p'}[0,T]$ by
$$
g(u) = |{_{0}}D_{t}^{\alpha}u|^{p-2}{_{0}}D_{t}^{\alpha}u,
$$
for $u\in E_{0}^{\alpha ,p}$. Let us prove that $g$ is continuous.  

\noindent 
{\bf Case $p \in (2,\infty)$.} For $u,v \in E_{0}^{\alpha,p}$, by (\ref{MR03}) and H\"older inequality we have:
\begin{eqnarray}\label{MR05}
\int_{0}^{T} |g(u) - g(v)|^{p'} dt &=& \int_{0}^{T} \left| |{_{0}}D_{t}^{\alpha}u|^{p-2}{_{0}}D_{t}^{\alpha}u - |{_{0}}D_{t}^{\alpha}v|^{p-2}{_{0}}D_{t}^{\alpha}v   \right|dt\nonumber\\
&\leq & \beta \int_{0}^{T} |{_{0}}D_{t}^{\alpha}u - {_{0}}D_{t}^{\alpha}v|^{p'} \left( |{_{0}}D_{t}^{\alpha}u| + |{_{0}}D_{t}^{\alpha}v| \right)^{p'(p-2)}dt\nonumber\\
&\leq& \beta \left( \int_{0}^{T} |{_{0}}D_{t}^{\alpha}u - {_{0}}D_{t}^{\alpha}v|^{p} dt \right)^{p'/p} \left( \int_{0}^{T} \left[ |{_{0}}D_{t}^{\alpha}u| + |{_{0}}D_{t}^{\alpha}v|  \right]^{p}dt  \right)^{\frac{p'(p-2)}{p}}\nonumber\\
& = & \beta \|{_{0}}D_{t}^{\alpha}u - {_{0}}D_{t}^{\alpha}v\|_{L^p}^{p'} \||{_{0}}D_{t}^{\alpha}u| + |{_{0}}D_{t}^{\alpha}v|\|_{L^p}^{p'(p-2)}\nonumber\\
&\leq& \overline{C} \|u-v\|_{\alpha ,p}^{p'} \left(\|u\|_{\alpha ,p} + \|v\|_{\alpha,p}  \right)^{p'(p-2)} 
\end{eqnarray}
with $\overline{C}$ constant independent of $u$ and $v$.

\noindent 
{\bf Case $p\in (1,2]$.} For $u,v\in E_{0}^{\alpha,p}$, by (\ref{MR04}) it follows
\begin{eqnarray}\label{MR06}
\int_{0}^{T} |g(u) - g(v)|^{p'}dt & = & \int_{0}^{T} \left| |{_{0}}D_{t}^{\alpha}u|^{p-2} {_{0}}D_{t}^{\alpha}u - |{_{0}}D_{t}^{\alpha}v|^{p-2}{_{0}}D_{t}^{\alpha}v  \right|^{p'}dt\nonumber\\
&\leq& \beta \int_{0}^{T} |{_{0}}D_{t}^{\alpha}u - {_{0}}D_{t}^{\alpha}v|^{p'(p-1)}dt\nonumber\\
&\leq& \overline{C}_1\|u-v\|_{\alpha ,p}^{p-1}
\end{eqnarray}
with $\overline{C}_1$ constant independent of $u$ and $v$. From (\ref{MR05}) and (\ref{MR06}) the continuity of $g$ is obvious.

On the other hand, we claim that
\begin{equation}\label{MR07}
\|J'(u) - J'(v)\|_{(E_{0}^{\alpha , p})^{*}} \leq K\|g(u) - g(v)\|_{L^{p'}}
\end{equation}
with $K>0$ constant independent of $u,v \in E_{0}^{\alpha,p}$. Indeed, by the H\"older inequality we have:
\begin{eqnarray*}
\left| \langle J'(u) - J'(v), \varphi \rangle \right| &\leq& \int_{0}^{T} |g(u) - g(v)||{_{0}}D_{t}^{\alpha} \varphi|dt\\
&\leq& \left( \int_{0}^{T} |g(u) - g(v)|^{p'}  \right)^{\frac{1}{p'}} \left( \int_{0}^{T} |{_{0}}D_{t}^{\alpha}\varphi|^pdt \right)^{\frac{1}{p}}\\
&\leq& K \|g(u) - g(v)\|_{L^{p'}}\|\varphi\|_{\alpha ,p}
\end{eqnarray*} 
for $u,v,\varphi \in E_{0}^{\alpha ,p}$, proving (\ref{MR07}).

Now, by the continuity of $g$ and (\ref{MR07}), the conclusion of the Lemma follows in a standard way. $\Box$ 

Combining Lemma \ref{MR1lem} and Lemma \ref{MR2lem}, we get that $I\in C^1(E_{0}^{\alpha ,p}, \mathbb{R})$ and
$$
\langle I'(u),v \rangle = \int_{0}^{T} |{_{0}}D_{t}^{\alpha}u|^{p-2}{_{0}}D_{t}^{\alpha}u {_{0}}D_{t}^{\alpha}vdt - \int_{0}^{T} f(t,u)vdt,
$$
for all $u,v\in E_{0}^{\alpha ,p}$.

\begin{Lem}\label{PT1lem}
Suppose that $f$ satisfies $(f_{1}) - (f_{3})$. Then there exist $\rho >0$ and $\beta >0$ such that
$$
I(u) \geq \alpha >0,
$$
for any $u\in E_{0}^{\alpha ,p}$ with $\|u\|_{\alpha ,p} = \rho$.
\end{Lem}

\noindent 
{\bf Proof.} By assumptions $(f_1)$ and $(f_3)$, for any $\epsilon >0$, there exists $C_{\epsilon} >0$ such that for any $\xi \in \mathbb{R}$ and a.e. $t\in [0,T]$, we have
\begin{equation}\label{MR08}
|f(t,\xi)| \leq p\epsilon|\xi|^{p-1} + qC_{\epsilon} |\xi|^{q-1}.
\end{equation} 
It follows from (\ref{MR08}) that
\begin{equation}\label{MR09}
|F(t,\xi)| \leq \epsilon |\xi|^{p} + C_\epsilon |\xi|^{q}.
\end{equation}
Let $u\in E_{0}^{\alpha ,p}$. By (\ref{MR09}), Proposition \ref{FC-FEprop3} and Remark \ref{embb}, we obtain
\begin{eqnarray}\label{MR10}
I(u) & = & \frac{1}{p}\int_{0}^{T}|{_{0}}D_{t}^{\alpha}u(t)|^pdt - \int_{0}^{T}F(t,u(t))dt\nonumber\\
&\geq& \frac{1}{p} \int_{0}^{T} |{_{0}}D_{t}^{\alpha}u(t)|^{p}dt - \epsilon \int_{0}^{T}|u(t)|^pdt - C_{\epsilon}\int_{0}^{T}|u(t)|^{q}dt\nonumber\\
&\geq& \frac{1}{p}\|u\|_{\alpha ,p}^{p} - \frac{\epsilon T^{\alpha}}{\Gamma (\alpha +1)}\|u\|_{\alpha ,p}^{p} + C_\epsilon \mathcal{K}\|u\|_{\alpha ,p}^{q}
\end{eqnarray}
where 
$$\mathcal{K} = \frac{T^{\alpha q + 1 - \frac{q}{p}}}{(\Gamma (\alpha)[(\alpha - 1)q + 1]^{1/q})^{q-p}\Gamma (\alpha + 1)^{p}}.$$
Choosing $\epsilon = \frac{\Gamma (\alpha +1)}{2pT^{\alpha}}$, by (\ref{MR10}), we have
$$
I(u) \geq \frac{1}{2p}\|u\|_{\alpha ,p}^{p} - C\|u\|_{\alpha,p}^{q} \geq \|u\|_{\alpha,p}^{p}\left( \frac{1}{2p} - C\|u\|_{\alpha,p}^{q-p}  \right),
$$ 
where $C$ is a constant only depending on $\alpha, p, T$. Now, let $\|u\|_{\alpha,p} = \rho >0$. Since $q>p$, we can choose $\rho$ sufficiently small such that
$$
\frac{1}{2p} - C\rho^{q-p} >0,
$$
so that
$$
I(u) \geq \rho^{p}\left(\frac{1}{2p} - C\rho^{q-p}  \right) =:\beta >0.
$$
Thus, the Lemma is proved. $\Box$
\begin{Lem}\label{PT2lem}
Suppose that $f$ satisfies $(f_{1}) - (f_{3})$. Then there exists $e\in C_{0}^{\infty}[0,T]$ such that $\|e\|_{\alpha ,p} \geq \rho$ and $I(e) < \beta$, where $\rho$ and $\beta$ are given in Lemma \ref{PT1lem}. 
\end{Lem}

\noindent
{\bf Proof.} From assumption ($f_2$) it follows that
\begin{equation}\label{MR11}
F(t, \xi) \geq r^{-\mu} \min\{F(t,r), F(t,-r)\}|\xi|^{\mu}
\end{equation}
for all $|\xi| >r$ and a.e. $t\in [0,T]$. Thus, by (\ref{MR11}) and $F(t,\xi) \leq \max_{|\xi| \leq r}F(t,\xi)$ for all $|\xi| \leq r$, we obtain
\begin{equation}\label{MR12}
F(t,\xi) \geq r^{-\mu} \min\{F(t,r), F(t,-r)\}|\xi|^{\mu} - \max_{|\xi|\leq r} F(t,\xi) - \min \{F(t,r), F(t,-r)\},
\end{equation}
for any $\xi \in \mathbb{R}$ and a.e. $t\in [0,T]$.
Since $C_{0}^{\infty}[0,T] \subset E_{0}^{\alpha,p}$, we can fix $u_{0} \in C_{0}^{\infty}[0,T]$ such that $\|u_0\|_{\alpha ,p} = 1$. Now, let $\sigma \geq 1$, by (\ref{MR12}), we have 
\begin{eqnarray*}
I(\sigma u_0) & = & \frac{\sigma^{p}}{p}\|u_0\|_{\alpha ,p}^{p} - \int_{0}^{T}F(t,\sigma u_0(t))dt\\
&\leq& \frac{\sigma^p}{p} - r^{-\mu}\sigma^{\mu} \int_{0}^{T}\min \{F(t,r), F(t,-r)\}|u_0(t)|^{\mu}dt\\
&&+ \int_{0}^{T} \max_{|\xi|\leq r}F(t,\xi) + \min \{F(t,r), F(t,-r)\}dt.
\end{eqnarray*}
From assumption $(f_1)$ and $(f_2)$, we get that $0< F(t,\xi) \leq C(|r| + |r|^{q})$ for $|\xi| \leq r$ a.e.$t\in [0,T]$. Thus,
$$
0< \min \{F(t, r), F(t,-r)\} < C(|r| + |r|^q) \;\;\mbox{a.e.}\; t\in [0,T].
$$
Since $\mu >p$, passing to the limit as $t\to \infty$, we obtain that $I(tu_0) \to -\infty$. Thus, the assertion follows by taking $e = Tu_0$ with $T$ sufficiently large. $\Box$

\begin{Lem}\label{PT3lem}
Suppose that $f$ satisfies $(f_{1}) - (f_{3})$. Then the functional $I$ satisfies (PS) condition.
\end{Lem}

\noindent
{\bf Proof.} For any sequence $\{u_n\}\subset E_{0}^{\alpha,p}$ such that $I(u_n)$ is bounded and $I'(u_n) \to 0$ as $n\to \infty$, there exists $M>0$ such that
$$
|\langle I'(u_n), u_n \rangle| \leq M\|u_n\|_{\alpha,p}\;\;\mbox{and}\;\;|I(u_n)| \leq M.
$$ 
For each $n\in \mathbb{N}$, we denote
$$
\Omega_{n} = \{t\in [0,T]|\;|u_n(t)| \geq r\},\;\Omega'_n = [0,T]\setminus \Omega_n.
$$
We have
\begin{equation}\label{MR13}
\frac{1}{p}\|u_n\|_{\alpha,p}^{p} - \left( \int_{\Omega_n}F(t,u_n) + \int_{\Omega'_n} F(t,u_n) \right) \leq M.
\end{equation}
We proceed with obtaining estimations independent of $n$ for the integrals in (\ref{MR13}). Let $n\in \mathbb{N}$ be arbitrary chosen. From assumption ($f_1$), we have
\begin{equation}\label{MR14}
|F(t,\xi)| \leq 2C(|\xi|^{q} + 1).
\end{equation}
If $t\in \Omega'_n$, then $|u_n(t)| <r$ and by (\ref{MR14}), it follows
$$
F(t,u_n) \leq 2C(|u_n|^{q} + 1)\leq 2C(r^q + 1)
$$
and hence
\begin{equation}\label{MR15}
\int_{\Omega'_n} F(t,u_n)dt \leq 2CTr^q  + T = K_1.
\end{equation}
If $t\in \Omega_n$, then $|u_n(t)| \geq r$ and by ($f_2$) it holds
$$
F(t,u_n) \leq \frac{1}{\mu}f(t,u_n(t))u_n(t)
$$
which gives
\begin{equation}\label{MR16}
\int_{\Omega_n} F(t,u_n)dt \leq \int_{\Omega_n} \frac{1}{\mu}f(t,u_n(t))u_n(t)dt = \frac{1}{\mu} \left( \int_{0}^{T} f(t,u_n)u_ndt - \int_{\Omega'_n} f(t,u_n)u_ndt \right)
\end{equation}
By ($f_1$), we deduce
\begin{eqnarray*}
\left| \int_{\Omega'_n}f(t,u_n)u_ndt \right| &\leq& \int_{\Omega'_n} C(|u_n| + |u_n|^q)dt\\
&\leq& CTr + CTr^q = K_2,
\end{eqnarray*}
which yields 
\begin{equation}\label{MR17}
-\frac{1}{\mu}\int_{\Omega'_n}f(t,u_n)u_ndt \leq \frac{K_2}{\mu}.
\end{equation}
Finally, by (\ref{MR13}), (\ref{MR15}), (\ref{MR16}) and (\ref{MR17}) we obtain 
\begin{eqnarray}\label{MR18}
\frac{1}{p}\|u_n\|_{\alpha,p}^{p} - \frac{1}{\mu}\int_{0}^{T}f(t,u_n)u_ndt &\leq& M + K_1 + \frac{K_2}{\mu} = K,\nonumber\\
\frac{1}{p}\|u_n\|_{\alpha,p}^{p} - \frac{1}{\mu}\langle H'(u_n), u_n\rangle &\leq& K 
\end{eqnarray} 
On the other hand, since $|\langle I'(u_n), u_n \rangle| \leq M \|u_n\|_{\alpha ,p}$ for $n\geq n_0$. Consequently, for all $n\geq n_0$, we have
$$
|\|u_n\|_{\alpha,p}^{p} - \langle H'(u_n), u_n\rangle| \leq M\|u_n\|_{\alpha,p}
$$
which gives
\begin{equation}\label{MR19}
-\frac{1}{\mu}\|u_n\|_{\alpha ,p}^{p} - \frac{M}{\mu} \|u_n\|_{\alpha,p} \leq -\frac{1}{\mu} \langle H'(u_n), u_n \rangle.
\end{equation}
Now, from (\ref{MR18}) and (\ref{MR19}) it results
$$
\left(\frac{1}{p} - \frac{1}{\mu}  \right)\|u_n\|_{\alpha ,p}^{p} - \frac{M}{\mu}\|u_n\|_{\alpha,p} \leq K
$$
and taking into account that $\mu >p$, we conclude that $\{u_n\}$ is bounded. Since $E_{0}^{\alpha,p}$ is a reflexive Banach space, up to a subsequence, still denoted by $\{u_n\}$ such that $u_n \rightharpoonup u$ in $E_{0}^{\alpha,p}$. Then $\langle I'(u_n), u_n-u \rangle \to 0$. Thus, we obtain
\begin{eqnarray}\label{MR20}
\langle I'(u_n), u_n - u \rangle & = & \int_{0}^{T} |{_{0}}D_{t}^{\alpha}u_n|^{p-2}{_{0}}D_{t}^{\alpha}u_n({_{0}}D_{t}^{\alpha}u_n-{_{0}}D_{t}^{\alpha}u)dt - \int_{0}^{T}f(t,u_n)(u_n-u)dt\nonumber\\
&\to& 0
\end{eqnarray}
as $n\to \infty$. Moreover, by Proposition \ref{FC-FEprop4}, 
\begin{equation}\label{MR21}
u_n \to u\;\mbox{strongly in}\; C[0,T].
 \end{equation}
 From (\ref{MR21}), $\{u_n\}$ is bounded in $C[0,T]$, the by assumption ($f_1$), we have
 \begin{eqnarray*}
 \left| \int_{0}^{T} f(t,u_n)(u_n - u)dt \right| &\leq& \int_{0}^{T}|f(t,u_n)||u_n - u|dt\\
 &\leq& C\int_{0}^{T} |u_n - u|dt + C\int_{0}^{T} |u_n|^{q-1}|u_n - u|dt\\
 &\leq& CT\|u_n - u\|_{\infty} + CT\|u_n\|_{\infty}^{q-1}\|u_n - u\|_{\infty}.
 \end{eqnarray*}
 This combined with (\ref{MR21}) follows
 $$
\lim_{n\to \infty} \int_{0}^{T}f(t,u_n)(u_n-u)dt = 0,
 $$
 hence one has
 \begin{equation}\label{MR22}
 \int_{0}^{T} |{_{0}}D_{t}^{\alpha}u_n|^{p-2}{_{0}}D_{t}^{\alpha}u_n({_{0}}D_{t}^{\alpha}u_n - {_{0}}D_{t}^{\alpha}u)dt \to 0\;\mbox{as}\;n\to \infty.
 \end{equation} 
 Using the standard inequality given by
 \begin{eqnarray*}
&& (|z|^{p-2}z - |y|^{p-2}y)(z-y) \geq C_{p}|z-y|^{p}\;\;\mbox{if}\;\;p\geq 2\\
&&(|z|^{p-2}z - |y|^{p-2}y)(z-y) \geq \tilde{C}_p\frac{|z-y|^2}{(|z| + |y|)^{2-p}}\;\;\mbox{if} \;\;1<p<2.
 \end{eqnarray*}  
 (see \cite{IPe}). From which we obtain for $p>2$
 \begin{eqnarray}\label{MR23}
 \int_{0}^{T} |{_{0}}D_{t}^{\alpha}u_n - {_{0}}D_{t}^{\alpha}u|^pdt &\leq& \frac{1}{C_p} \int_{0}^{T} \left[ |{_{0}}D_{t}^{\alpha}u_n|^{p-2}{_{0}}D_{t}^{\alpha}u_n - |{_{0}}D_{t}^{\alpha}u|^{p-2}{_{0}}D_{t}^{\alpha}u \right] ({_{0}}D_{t}^{\alpha}u_n - {_{0}}D_{t}^{\alpha}u)dt\nonumber\\
  &&\to 0,
 \end{eqnarray}
 as $n\to \infty$. For $1<p<2$, by reverse H\"older inequality, we have
 \begin{eqnarray}\label{MR24}
\int_{0}^{T} |{_{0}}D_{t}^{\alpha}u_n - {_{0}}D_{t}^{\alpha}u|^pdt & \leq & \tilde{C}_{p}^{-\frac{p}{2}} \left( \int_{0}^{T}(|{_{0}}D_{t}^{\alpha}u_n| + |{_{0}}D_{t}^{\alpha}u|)^pdt \right)^{\frac{2-p}{2}}\nonumber\\
&\times&\!\!\!\!\! \left(\int_{0}^{T}[|{_{0}}D_{t}^{\alpha}u_n|^{p-2}{_{0}}D_{t}^{\alpha}u_n - |{_{0}}D_{t}^{\alpha}u|^{p-2}{_{0}}D_{t}^{\alpha}u]({_{0}}D_{t}^{\alpha}u_n - {_{0}}D_{t}^{\alpha}u)dt\right)^{p/2}\nonumber\\
&\leq&\!\!\!\! \overline{C} \left(\int_{0}^{T}[|{_{0}}D_{t}^{\alpha}u_n|^{p-2}{_{0}}D_{t}^{\alpha}u_n - |{_{0}}D_{t}^{\alpha}u|^{p-2}{_{0}}D_{t}^{\alpha}u][{_{0}}D_{t}^{\alpha}u_n - {_{0}}D_{t}^{\alpha}u]dt  \right)^{p/2}\nonumber\\
&\to &0,
 \end{eqnarray}
 as $n\to \infty$. Combing (\ref{MR23}) with (\ref{MR24}), we get that $u_n \to u$ strongly in $E_{0}^{\alpha,p}$ as $n \to \infty$. Therefore, $I$ satisfies (PS) condition. $\Box$
 
 \noindent 
 {\bf Proof of Theorem \ref{main}.}
 Since Lemma \ref{PT1lem} - Lemma \ref{PT3lem} hold, the Mountain pass Theorem (see \cite{PR}) gives that there exists a critical point $u\in E_{0}^{\alpha ,p}$ of $I$. Moreover,
 $$
 I(u) \geq \beta > 0 = I(0).
 $$ 
 Thus, $u\neq 0$. $\Box$
 





\end{document}